\documentclass[11pt]{article}
\usepackage{amsmath}
\usepackage{amsthm}
\usepackage{amssymb}

\newtheorem{theo}{Theorem}[section]

\newtheorem{prop}[theo]{Proposition}

\newtheorem{quest}[theo]{Question}

\newcommand{\FF}{{\cal F}}
\newcommand{\GG}{{\cal G}}
\newcommand{\HH}{{\cal H}}

\newcommand{\KK}{{\cal K}}

\begin{document}
\date{}

\title{Graph-Codes}

\author{Noga Alon
\thanks{Princeton University,
Princeton, NJ, USA 
and
Tel Aviv University, Tel Aviv,
Israel.
Email: {\tt nalon@math.princeton.edu}.
Research supported in part by
NSF grant DMS-2154082 and by USA-Israel BSF grant 2018267.}
}

\maketitle
\begin{abstract}
The symmetric difference of two graphs $G_1,G_2$ on the same set of
vertices $[n]=\{1,2, \ldots ,n\}$ is the graph on $[n]$ whose set
of edges are all edges that belong to exactly one of the two graphs
$G_1,G_2$. Let $H$ be a fixed graph with an even (positive) number
of edges, and let $D_H(n)$ denote the maximum possible cardinality
of a family of graphs on $[n]$ containing no two members whose 
symmetric difference is a copy of $H$. Is it true that 
$D_H(n)=o(2^{n \choose 2})$ for any such $H$? We discuss this
problem, compute the value of $D_H(n)$ up to a constant factor for
stars and matchings, and discuss several variants of the problem
including ones that have been considered in earlier work.
\end{abstract}

\section{Introduction}

\subsection{The problem}

The {\em symmetric difference} of two graph $G_1=(V,E_1)$ and $G_2=(V,E_2)$
on the same set of vertices $V$ is the graph $(V,E_1\oplus E_2)$
where $E_1 \oplus E_2$ is the
symmetric difference between $E_1$
and $E_2$, that is, the set of all edges that belong to exactly one
of the two graphs. Put $V=[n]=\{1,2, \ldots ,n\}$ and let
$\HH$ be a family of graphs on the set of vertices $[n]$ 
which is closed under
isomorphism.
A collection of graphs $\FF$ on $[n]$ is called 
an {\em $\HH$-(graph)-code} 
if it contains no two members whose symmetric difference is a graph
in $\HH$. For the special case that $\HH$ contains all copies 
of a single
graph $H$ on $[n]$ this is called an $H$-code. 
Here we are interested
in the maximum possible cardinality of such codes for various families 
$\HH$. Let $D_{\HH}(n)$ denote this maximum, and let
$$
d_{\HH}(n)=\frac{D_{\HH}(n)}{2^{n \choose 2}}
$$
denote the maximum possible fraction of the total number of graphs
on $[n]$ in an $\HH$-code. If $\HH$ consists of all graphs
isomorphic to one graph $H$, we denote $d_{\HH}(n)$ by $d_{H}(n)$.
Note that if $\HH$ consists of all graphs with less than $d$ edges,
then
$D_{\HH}(n)$ is simply the maximum possible cardinality of a binary
code of length ${n \choose 2}$ and minimum distance at least $d$.
This motivates the terminology ``graph-codes'' used here.

The case 
$\HH=\KK$ where $\KK$ is the family
of all cliques is of particular interest. This
case is motivated by a conjecture of Gowers raised in his blog post
\cite{Go} in 2009 and is discussed briefly in the comments of that blog.
If $\HH$ consists of all graphs  with independence number at most
$2$, then $d_{\HH}(n) \geq 1/8$ for all $n \geq 3$, as shown by the
family of all graphs on $[n]$ containing a triangle on the set of
vertices $\{1,2,3\}$. An interesting result of Ellis, Filmus and 
Friedgut \cite{EFF}, settling a conjecture of Simonovits and
S\'os, asserts that this is tight for all $n \geq 3$. The
corresponding result, that $d_{\HH'}(n)=1/2^6$ for all 
$n \geq 4$, where $\HH'$ is the family of  all graphs with
independence number at most $3$, is proved in \cite{BZ}.  
A more systematic study of the parameters 
$D_{\HH}(n)$ and $d_{\HH}(n)$ for various families of graphs
$\HH$  appears in the recent paper \cite{AGKMS}. The families
$\HH$ considered in this work include the family of all
disconnected graphs, the family of all graphs that are not
$2$-connected, the family of all non-Hamiltonian graphs and the
family of all graphs that contain or do not contain a spanning star. 
Additional families studied are all graphs that contain an induced
or non-induced copy of a
fixed graph $T$, or all graphs that do not contain such a subgraph.

In this note we focus on the case that $\HH$  consists of a single
graph $H$ and the case that $\HH$ is the family of all cliques, or
all cliques up to a prescribed size. 
Note that trivially, if every member of $\HH$ has an odd number of edges
then $d_{\HH}(n) \geq \frac{1}{2}$ as the family of all graphs on
$[n]$ with an even number of edges forms an $\HH$-code.

This suggests the following intriguing question.
\begin{quest}
\label{q15}
Let $\HH$ be a  family of graphs closed under isomorphism.
Is it true that $d_{\HH}(n)$ tends to $0$ as $n$ tends to infinity
if and only if $\HH$ contains a graph with an even number of edges ?
Equivalently: is it true that for any fixed graph $H$ with an even
number of edges, $d_H(n)$ tends to $0$ as $n$ tends to infinity ?
\end{quest}

We also study the linear
variant of these problems, where the $\HH$-codes considered 
are restricted to linear subspaces, that is, to families of graphs
on $[n]$ closed under symmetric difference.

\subsection{Results}

Recall that $\KK$ is the family of all cliques.
Let $\KK(r)$ denote the set of all cliques on at most $r$ vertices.
Let $K_{1,t}$ denote the star with $t$ edges and let 
$M_t$ denote the matching of $t$ edges.

\begin{theo}
\label{t11}
For every positive integer $k$,
$$
d_{K_{1,2k}}(n)=\Theta_k(1/n^k)~~\mbox{and}~~~
d_{M_{1,2k}}(n)=\Theta_k(1/n^k).
$$
\end{theo}

\begin{prop}
\label{p12}
For every integer $r \geq 1$,
$$d_{\KK(4r+3)}(n)  \geq \Omega(\frac{1}{n^r}).$$
\end{prop}

\begin{prop}
\label{p13}
For the family $\KK$ of all cliques,
$d_{\KK}(n) \geq \frac{1}{2^{[n/2]}}$.
\end{prop}

\begin{prop}
\label{p14}
Let $H$ be a fixed graph obtained from two copies of a graph $H'$ 
by identifying the vertices of an independent set of $H'$. Then
$$
d_{H}(n) \leq \frac{|V(H)|}{n} ~~\mbox{for all}~~
n \geq |V(H)|.
$$ 
In particular, 
$d_H(n)$ tends to $0$ as $n$ tends to infinity.
\end{prop}
\vspace{0.5cm}

\noindent
{\bf Remark:}\, all lower bounds are proved by exhibiting  
proper colorings of the relevant Cayley graphs, and in all cases 
the constructed family is an affine space over $Z_2$. Using 
a simple Ramsey-theoretic argument  it is
not difficult to show that for an affine space the maximum
possible cardinality obtained is at most a fraction 
$O(\log \log n/\log n)$ of all graphs on $n$ vertices whenever
the defining family contains a fixed graph with an even 
number of edges.
\vspace{0.2cm}

\noindent
Since all lower bounds are obtained by what may be called linear 
graph-codes one can study this separately, as done for standard 
error correcting codes. For the family of all cliques $\KK$ we get here
an exact result (strengthening the assertion of Proposition \ref{p13}).
\begin{theo}
\label{t16}
For any $n \geq 2$, the
minimum possible co-dimension of a linear space of graphs
on $n$ vertices that contains no member of $\KK$ is
exactly $[n/2]$.
\end{theo}

\section{Proofs}

\subsection{Upper bounds}

For a family of graphs $\HH$  and an integer $n$, the Cayley
graph $C(n,\HH)$ is the graph whose vertices are all graphs
on the $n$ vertices $[n]$, where two are adjacent iff their symmetric
difference is a member of $\HH$. This is clearly a Cayley graph 
over the elementary abelian $2$-group $Z_2^{N}$ with $N={n \choose 2}$.
The function $D_{\HH}(n)$ is just the independence number of this graph,
$d_{\HH}(n)$ is the so called independence ratio. Since the graph
$C(n,\HH)$ is vertex transitive, its independence ratio is exactly the
reciprocal of its fractional chromatic number. In order to prove an 
upper bound of $\alpha$ for its independence ratio it suffices 
to exhibit a set $S$ of vertices that contains
no independent set of size larger than $\alpha |S|$. This applies
also to weighted sets of vertices, but we will not use weights here.
\vspace{0.5cm}

\noindent
{\bf Proof of Proposition \ref{p14}:}\, 
Let $a+b$ denote the number of vertices of $H'$ where $b$ is
the size of its independent set so that $H$ is obtained from two
copies of $H'$ by identifying the vertices in this independent set.
Thus the number of vertices of $H$ is $2a+b$. Consider the following
set of $m=\lfloor (n-b)/a \rfloor$ copies of $H'$ on subsets of the
vertex set $[n]$. All of them contain the same independent
set on the vertices $\{n-b+1,n-b+2, \ldots ,n\}$, and the additional
vertices of copy number $i$ are the vertices $\{(i-1)a+1,(i-1)a+2, 
\ldots ,ia\}$, where $1 \leq i \leq m$. Each of these copies can be
viewed as a vertex of the Cayley graph $C=C(n,\{H\})$. Since 
the symmetric difference of every pair
of such copies forms a copy of $H$, this set forms a clique of size
$m$ in $C$, implying that 
$d_{H}(n) \leq \frac{1}{m} \leq |V(H)|/n$. \hfill $\Box$
\vspace{0.5cm}

\noindent
The proofs of Theorem \ref{t11} for stars and for matchings 
are very similar. We describe the proof for stars and briefly
mention the modification needed for matchings. 
The upper bound in Theorem \ref{t11} for the star $K_{1,2}$ is a special
case of the result above (with $H'$ being a single edge).
The upper bound for any prime $k$ can be
proved using the following result of Frankl and Wilson.
\begin{theo}[\cite{FW}]
\label{tfw}
Let $p$ be a prime, and let $a_0,a_1, \ldots ,a_r$ be distinct
residue classes modulo $p$.
Let $\FF$ be a family of subsets of $[n]$  and suppose that
$|F| \equiv a_0 \bmod p$ for all $F \in \FF$ and  that for every
two distinct $F_1,F_2 \in \FF$,
$|F_1 \cap F_2 | \equiv a_i \bmod p$ for some $1 \leq i \leq r$.  
Then $|\FF| \leq \sum_{i=0}^r {n \choose i}$.
\end{theo}

Suppose $k$ is a prime, $n \geq 2k$ and consider the 
family $\GG$ of all stars $K_{1,2k-1}$ with center $1$
and $2k-1$ leaves  among the vertices $\{2,3, \ldots ,n\}$.
Thus $|\GG|={{n-1} \choose {2k-1}}$. If two such stars share exactly
$k-1$ common leaves then their symmetric difference is a copy of 
$K_{1,2k}$. A subset of $\GG$ which is independent in the Cayley
graph $C(n,K_{1,2k})$ corresponds to a collection of subsets of the
set $\{2,3,\ldots ,n\}$, each of size $2k-1$, where the
intersection of no two of these subsets is of cardinality $k-1$.
Therefore, each of these sets is of cardinality 
$-1$ modulo $k$ and no intersection is of cardinality $-1$ modulo
$k$. By the Frankl-Wilson Theorem (Theorem \ref{tfw})
the cardinality of such a family is at most
$\sum_{i=0}^{k-1}{{n-1} \choose i}$. Therefore, for every prime
$k$,
$$
d_{K_{1,2k}}(n) \leq 
\frac{{\sum_{i=0}^{k-1}{{n-1} \choose i}}}
{{{n-1} \choose {2k-1}}}
\leq O_k(\frac{1}{n^k}).
$$

In order to prove the upper bound for all $k$
we need the following result of Frankl and
F\"{u}redi.
\begin{theo}[\cite{FF}]
\label{t21}
For every fixed positive integers $\ell > \ell_1+\ell_2$ there exist
$n_0=n_0(\ell)$ and $d_{\ell}>0$ so that for all $n>n_0$,
if $\FF$ is a family of $\ell$-subsets of
$[n]$ in which the intersection of each pair of distinct
members is of cardinality either at least $\ell-\ell_1$ or strictly
smaller than $\ell_2$, then
$$
|\FF| \leq d_{\ell} \cdot n^{\max\{\ell_1, \ell_2\}}.
$$
\end{theo}
\vspace{0.5cm}

\noindent
{\bf Proof of Theorem \ref{t11}, upper bound:}\, 
The proof for stars is essentially identical to the one described above
for prime $k$, using Theorem \ref{t21} instead of Theorem
\ref{tfw}. 
Let $\GG$ be the family of all stars $K_{1,2k-1}$ with center $1$
and $2k-1$ leaves  among the vertices $\{2,3, \ldots ,n\}$.
Thus $|\GG|={{n-1} \choose {2k-1}}$. If two such stars share exactly
$k-1$ common leaves then their symmetric difference is a copy of $K_{1,2k}$.
Therefore, by Theorem \ref{t21} above with $\ell=2k-1, \ell_1=\ell_2=k-1$,
the maximum cardinality of 
a subset of $\GG$ which is independent in the Cayley graph
$C(n,K_{1,2k})$ is at most some $c_k (n-1)^{k-1}$ 
for all sufficiently large
$n$. This supplies the required upper bound 
$$
\frac{c_k (n-1)^k}{|\GG|} \leq O_k(\frac{1}{n^k}),
$$
for $d_{K_{1,2k}}(n)$.  The proof for matchings is similar,
starting with the family of all subsets of cardinality $2k-1$
of a fixed matching of cardinality $\lfloor n/2 \rfloor.$ 
The symmetric difference of any two matchings that share exactly 
$k-1$ common edges is a copy of $M_{2k}$. Thus the proof can proceed
exactly as in the case of stars.
\hfill  $\Box$

\subsection{Lower bounds}

In order to lower bound the independence number of a Cayley graph
$C=C(n,\HH)$ it suffices to upper bound its chromatic number. One way
to do so is to
assign to each edge $e$ of the complete graph on $[n]$ a vector 
$v_e \in Z_2^r$ for some $r$, so that for every $H \in \HH$,
$\sum_{e \in E(H)} v_e \neq 0$, where the sum is computed in
$Z_2^r$. Given these vectors, we can assign to each graph
$G$ on $[n]$ the color $\sum_{e \in E(G)} v_e$ (computed, of course,
in $Z_2^r$). This is clearly a proper coloring of $C$ by at most
$2^r$ colors. Note that the matrix whose columns are the 
${n \choose 2}$ vectors $v_e$ is the analogue of the parity-check
matrix of a linear error correcting code in the traditional 
theory of codes, and the color defined above is the analogue of
the syndrome of a word, see, e.g., \cite{MS} for more information
about these basic notions.
\vspace{0.5cm}

\noindent
{\bf Proof of Theorem \ref{t11}, lower bound:}\,
For stars, it suffices to show that the chromatic number of the Cayley graph
$C=C(n,K_{1,2k})$ is at most $O(n^k)$. Let $s$ be the smallest
integer so that $2^{s}-1 \geq n$. As shown by the columns of the
parity check matrix of a BCH-code with designed distance 
$2k+1$ there is a collection $S$ of $2^s-1$ binary vectors  of
length $r=ks$ so that no sum of at most $2k$ of them (in $Z_2^{ks}$) is
the zero vector. Fix a proper edge coloring $c$ of $K_n$ by $n$ colors.
For each edge $e$ let $v_e$ be the vector number $c(e)$ in $S$.
This gives the desired lower bound for stars.  For matchings we use 
essentially the same construction, starting with a (non-proper)
edge coloring of $K_n$ by $n$ colors in which each color class
forms a star.
\hfill $\Box$ 
\vspace{0.5cm}

\noindent
{\bf Proof of Proposition \ref{p12}, lower bound:}\,
As in the previous proof, but the initial edge-coloring now is
defined by $c(ij)=i$ for all $i<j$ and the binary vectors selected are
taken from the columns of the parity check matrix of 
a code with designed distance $2r+2$. Let $U$ be the set
of vertices of a clique of size at least $2$ and at most
$4r+3$. Then $U$ contains at least $1$ and at most $2r+1$ vertices
$i$ for which there is an odd number of vertices of $U$ with
index strictly larger than $i$. Therefore the sum of vectors
corresponding to the edges of the clique on $U$ is equal to 
a sum of at most
$2r+1$ column vectors of the parity check matrix, 
which is nonzero. \hfill $\Box$
\vspace{0.5cm}

\noindent
{\bf Proof of Proposition \ref{p13}, lower bound:}\,
This follows from the construction in the proof of Theorem
\ref{t16} described in the next section.

\section{Linear graph-codes}

\noindent
{\bf Proof of Theorem \ref{t16}:}\,
The theorem is equivalent to the statement that for all $n \geq 2$
the minimum possible $r=r(n)$ so that there are graphs 
$G_1,\ldots ,G_r$ on the vertex set $[n]$ such that every clique
on a subset of cardinality at least $2$ of $[n]$ contains an odd
number of edges of at least one graph $G_i$, is $r=[n/2]$.
It clearly suffices to prove the upper bound for odd $n$ (that
imply the result for $n-1$) and the lower bound for even $n$
(implying the result for $n+1$). The upper bound is
described in what follows.  
Let $n \geq 3$ be odd.
Split the numbers $[n-1]
=\{1,2,\ldots ,n-1\}$ into the $(n-1)/2$ blocks 
$B_i=\{2i-1,2i\}$ for $1 \leq i \leq (n-1)/2$. Let
$G_i$ be the graph consisting of all edges of the 
$n-2i$ triangles with a common base $B_i$ 
on the vertices $B_i \cup \{j\}$ for
$2i <j \leq n$. Our family of graphs is the set of these
$(n-1)/2$ graphs $G_i$.
Let $K$
be an arbitrary clique on a subset $A$ of at least $2$ vertices in
$[n]$. 
If $A$ contains a full block $B_i$ for some $i$, then
it contains exactly $2x+1$ edges of $G_i$, where $x$ is the cardinality
of the intersection of $A$ with $\{2i+1,2i+2, \ldots ,n\}$.
As this is odd for all $x \geq 0$ we may assume 
that $A$ contains no block $B_i$. In this case,
let $j$ be the second largest element in $A$ (recall that
$|A| \geq 2$). Clearly $j \leq n-1$, hence it is contained in one
of the blocks $B_i$. But in this case $G_i$ contains exactly
one edge of the clique $K$, completing the proof of the upper bound.
Note that it is simple to give additional constructions with the same
properties as any set of graphs that spans the same subspace as the
graphs above will do. In particular, we can replace one of the graphs
$G_i$ by the complete graph $K_n$, which is the sum of all graphs
$G_i$. 

To prove the lower bound assume $n$ is even and let 
$G_1, \ldots G_{n/2-1}$ be a family of $n/2-1$ graphs on $[n]$.
We have to show that there is a clique on at least $2$ vertices containing
an even number of edges of each $G_i$. We show that in fact there is such
a clique on an even number of vertices.  To do so we apply the classical
theorem of Chevalley and Warning (cf., e.g., \cite{BS} or \cite{Sch}). 
Recall that  it
asserts that any system of polynomials with $n$ variables over
a finite field in which the number of variables exceeds the sum
of the degrees, which admits a solution, must admit another one
(in fact, the number of solutions is divisible by the characteristics).
Associate each vertex $i$ with a variable $x_i$ over $Z_2$ and consider the
following homogeneous system of polynomial equations over $Z_2$.
For each graph $G_s$ in our family,
$$
\sum_{ij \in E(G_s)}  x_ix_j =0.
$$
In addition, add the linear equation $\sum_{i=1}^n x_i=0$.

The sum of the degrees of the polynomials here is $2(n/2-1)+1=n-1$,
which is smaller than the number of variables. Since the system
is homogeneous it admits the trivial solution $x_i=0$ for all $i$.
Any other solution (which exists by the Chevalley Warning Theorem)
gives a clique on the set of vertices $\{i: x_i =1\}$ which is
nonempty, of even cardinality, and contains an even number of edges
(possibly zero) of each $G_i$. This establishes the lower bound
and completes the proof of Theorem \ref{t16}. \hfill $\Box$

\section{Concluding remarks and open problems}

\begin{itemize}
\item
Question \ref{q15}, which is equivalent to the problem of
deciding whether or not
for any fixed nonempty graph $H$ with an even number of edges
$d_{H}(n)$ tends to $0$ as $n$ tends to infinity, remains wide
open. 

An interesting special case is whether or not $d_{K_4}(n)=o(1)$.
It is also interesting to decide whether or not
$d_{K_4}(n)\geq \frac{1}{n^{o(1)}}.$ It is not difficult to show
that the latter can be deduced  from the existence 
of an edge coloring of $K_n$ by $n^{o(1)}$ colors with no copy of
$K_4$ in which every color appears an even number of times.
Indeed, such a coloring together with the columns of the  
parity check matrix of a BCH code with designed distance
$7$ supplies the lower bound above using the reasoning in the
proofs of some of the results here. I have just learned
from Zach Hunter and Dhruv Mubayi 
that such an edge coloring is described in
\cite{CH}, modifying the constructions in 
\cite{Mu}, \cite{CFLS}.  Therefore
$d_{K_4}(n)\geq \frac{1}{n^{o(1)}}.$
\item
Gowers conjectured in \cite{Go} that any family of a constant 
fraction of all graphs on $[n]$, where $n$ is sufficiently large,
contains two graphs $G_1,G_2$ such that $G_2$ is a subgraph of
$G_1$ and the symmetric difference of the two graphs (that is, the
set of all edges of $G_1$ that are not in $G_2$) forms a clique.
This is clearly stronger than the conjecture that
$d_{\KK}(n)$ tends to $0$ as $n$ tends to infinity, which is also
open. As explained in \cite{Go} the question of Gowers can be
viewed as the first unknown case  of a 
polynomial version of the density Hales-Jewett Theorem.
\item
As mentioned in the remark following the statement of Proposition
\ref{p14}, it is not difficult to show that for every graph $H$
with an even number of eges the maximum possible
cardinality of a {\em linear} family of graphs on $[n]$ in which no
symmetric difference is a copy of $H$, is $o(2^{n \choose 2}).$
As the proof applies Ramsey's Theorem, it provides very weak
bounds. It will be interesting to  establish tighter bounds for
the linear case. Theorem \ref{t16} provides an example of a tight
result of this form.
\item
The problem considered above can be extended to hypergraphs. More
generally, it can be extended to other versions of problems about
binary codes, where the coordinates of each codeword are indexed by 
the elements of some combinatorial structure, and the forbidden
symmetric differences correspond to a prescribed family of
substructures. Here is an example of a problem of this type. 
What is the maximum possible cardinality of
a collection of binary vectors whose coordinates are indexed
by the elements of the ordered set $[n]$,
where no symmetric difference of two distinct members of the
collection forms an interval of length which is a cube of an
integer? The corresponding Cayley graph here has $2^n$ vertices,
and it is triangle-free by Fermat's last
Theorem for cubes. Its independece number, which is 
the answer to the question above, is
$o(2^n)$. Indeed, this follows from the Furstenberg-S\'ark\"ozy 
Theorem and its extensions \cite{Sa}, by considering the maximum possible
cardinality of an independent set in the induced subgraph on the
set of all vertices that are characteristic vectors of an interval
$[i]=\{1, \ldots ,i\}$ for $0 \leq i \leq n$. 
\item
Some of the discussion here suggests the problem of determining or estimating
the smallest number of colors in an edge coloring of $K_n$ in which every
copy of a given graph $H$ (or every copy of any member of a prescribed
family $\HH$ of graphs) intersects at least one of the color classes by an odd
number of edges. This appears to be an interesting variant of classical
questions in Ramsey Theory and deserves further study.
\end{itemize}
\vspace{0.2cm}

\noindent
{\bf Acknowledgment} I would like to thank Zach Hunter and Dhruv Mubayi 
for helpful
comments and in particular for telling me about \cite{CH}.

\end{document}